\documentclass[11pt]{article}
\usepackage{yhmath,graphicx,amsmath,amsfonts,amssymb,amsthm,color,amscd,mathtools}
\usepackage[normalem]{ulem}

\theoremstyle{plain}
\newtheorem{theorem}{Theorem}
\newtheorem{acknowledgment}{Acknowledgment}

\newtheorem{lemma}{Lemma}

\newtheorem{problem}{Problem}
\newtheorem{openprob}{Open Problem}
\newtheorem{proposition}{Proposition}

\numberwithin{equation}{section}

\begin{document}

\title{\textbf{Stable intersection of middle}-$ \alpha$  \textbf{Cantor sets}}
\date{}
\author{\textbf{M. Pourbarat\footnote{E-mail: m-pourbarat@sbu.ac.ir}}\\
{\begin{small}\it\hspace*{-.2cm}Department of Mathematics,\end{small}}{\begin{small} \it Faculty of Mathematical Sciences,\end{small}}
{\begin{small}\it \end{small}}\\[-10pt]{\begin{small} \it Shahid Beheshti University, G.C.,  Evin, Tehran 19839,
 Iran
\end{small}}\\[-20pt]} \maketitle

\begin{abstract}
In the present paper, We introduce a pair of middle Cantor sets namely $(C_\alpha, C_\beta)$ having stable intersection, while the product of their thickness is smaller than one. Furthermore, the arithmetic difference $C_\alpha- \lambda C_\beta$ contains at least one interval for each nonzero number $\lambda$.
\end{abstract}

\section{Introduction}
Regular Cantor sets appear in dynamical systems when hyperbolic sets intersect stable and unstable manifolds of its points. The study of metrical and topological properties of intersection of regular Cantor sets emerges naturally in the theory of homoclinic bifurcations. Meanwhile, stable intersections between Cantor sets which come from stable and unstable foliations of a horseshoe, provide examples of open sets of nonhyperbolic diffeomorphisms  after  unfolding  a homoclinic tangency. In these cases, the open set of diffeomorphisms presenting persistent tangencies between the stable and unstable foliations of the horseshoe stably has positive lower density at the initial parameter of the bifurcation, in parametrized families (see [M1] and [PT]).

Regular Cantor sets appear in number theory too,  related to diophantine approximations. Many Cantor sets given by combinatorial conditions on the continued fraction of real numbers emerge in this situation. In the study of the classical Markov and Lagrange spectra related to them,  we usually deal with the arithmetic difference of two regular Cantor sets (see [M2]). On the other hand, the intersection of regular Cantor sets can be interpreted by their arithmetic difference as
$$K - K':= \big\lbrace x-y\mid ~~x \in K~,~y \in K'~\big\rbrace = \big\lbrace t\in \mathbb{R}\mid~~ K\cap (K'+t)\neq \emptyset  ~\big\rbrace.$$

The arithmetic difference of regular Cantor sets  can be so complicate even for the simplest possible family of pairs of Cantor sets (see [MO] and [P3]). For instance, conjecture of Palis remains still open in the affine case (see [P1], [P2] and [PT]).
Also, there exist regular Cantor sets $K$ and $K'$, such that $K-K'$ has positive Lebesgue measure, but does not contain any interval  (see [S]).
Having stable intersection of regular Cantor sets, clearly implies the existence of an interval contained in their arithmetic difference.

Before stating our main result, we pose some notations. A Cantor set
$K$ is regular or dynamically defined if:

{\rm i) there are disjoint compact intervals $K_{1},K_{2},\ldots,K_{r}$ such that
$K\subset K_{1}\cup\cdots\cup K_{r}$ and the boundary of each $K_{i}$ is
contained in $K$,

{\rm ii) there is a $C^{1+\epsilon}$ expanding map $\psi$
defined in a neighborhood of set $K_{1}\cup K_{2}\cup\cdots\cup K_{r}$ such that
$\psi(K_{i})$ is the convex hull of a finite union of some intervals $K_{j}$
satisfying:
\\$~~~~~~~~~$ii.1 For each $i$, $1\leq i\leq r$ and $n$ sufficiently
big, $\psi^{n}(K\cap K_{i})=K,$
\\$~~~~~~~~~$ii.2 $K=\bigcap_{n=0}^{\infty}\psi
^{-n}(K_{1}\cup K_{2}\cup\cdots\cup K_{r})$.\\
The set $\{K_{1},K_{2},\cdots,K_{r}\}$ is, by definition, a Markov partition for $K$, and
the set $D:=\bigcup_{i=1}^{r} K_{i}$ is the Markov domain of $K$.

A regular Cantor set is affine if $D\psi$ be constant on every interval $K_i$. The simplest kind of affine Cantor sets are middle-$\alpha$  Cantor sets that generalizes in the most natural way, the usual ternary Cantor set   which corresponds to p=3 in below definition.
\textbf{Definition.} Let $p>2$ and $\alpha:= 1-\frac{2}{p}$. Then
\emph{middle}-$\alpha$  Cantor set can be written as
$$C_\alpha :=~\Big\{~x\in \mathbb{R}~\mid~~x=(1-\frac{1}{p})\sum_{i=0}^\infty~\frac{a_i}{p^i}~~~,~~~a_i\in\{0,1\}~\Big\},$$
that is a regular Cantor set with the Markov partition $\{K_{1},K_{2}\}$ and expanding map
$$\phi(x):= ~\left\lbrace  \begin{array}
[c]{lcr}
~~px  ~~~~~~~~~~~~~~~  x \in K_1:=~[0,~\frac{1}{p}]\\
px-p+1~~~~~~  x \in K_2:=~ [1-\frac{1}{p},~1]
\end{array}. \right .$$

We say that the Cantor set $K$ is close on the topology $C^{1+\epsilon}$  to a Cantor set $\widetilde{K}$ with the Markov partition $\{\widetilde{K_1},\widetilde{K_2},\cdots,\widetilde{K_s}\}$
defined by expanding map $\widetilde{\psi}$ if and only if $r=s$, the extremes of $K_i$ are near the corresponding extremes
of $\widetilde{K_i},~ i=1,2,...,r$ and supposing $\psi\in C^{1+\epsilon}$ with Holder constant $C$, we must have $\widetilde{\psi}\in C^{1+\widetilde{\epsilon}}$ with Holder constant $\widetilde{C}$
such that $(\widetilde{C},\widetilde{\epsilon})$ is near $(C,\epsilon)$ and $\widetilde{\psi}$ is close to $\psi$
in the $C^1$ topology.

 \textbf{Definition.} Regular Cantor sets $K$ and $K'$ have \emph{stable intersection} if for any pair of regular Cantor sets $(\widetilde{K} ,\widetilde{K'})$ near $(K ,K')$, we have
$\widetilde{K}\cap \widetilde{K'} \neq\emptyset$.

Besides the Hausdorff dimension, there is another fractal invariant namely thickness introduced by Newhouse, that plays a relevant role in determining stable intersection of regular Cantor sets (see [N]).  Such thickness condition was generalized by Moreira in [M1] as follows:

\textbf{Definition.} Take $U$ be a bounded gap  of Cantor set $K$ and  $L_U$, $R_U$ be the intervals at its left  and its right, respectively, that separate it from the closest larger gaps.

\unitlength1.30mm \linethickness{1pt}
\begin{picture}(111.33,12.27)
\put(13.00,5.00){\makebox(0,0)[lc]{$K:$}}
\put(20.00,5.00){\line(1,0){5.00}}
\put(35.00,5.00){\line(1,0){10.00}}
\put(53.00,5.00){\line(1,0){6.00}}
\put(70.00,5.00){\line(1,0){5.00}}

\put(20.00,3.70){\makebox(0,0)[cb]{$|$}}
\put(25.30,3.70){\makebox(0,0)[cb]{$($}}
\put(34.70,3.70){\makebox(0,0)[cb]{$)$}}
\put(45.30,3.70){\makebox(0,0)[cb]{$($}}
\put(52.70,3.70){\makebox(0,0)[cb]{$)$}}
\put(59.30,3.70){\makebox(0,0)[cb]{$($}}
\put(69.70,3.70){\makebox(0,0)[cb]{$)$}}
\put(75.00,3.70){\makebox(0,0)[cb]{$|$}}

\put(40.00,7.00){\makebox(0,0)[cb]{$L_U$}}
\put(49.00,5.50){\makebox(0,0)[cb]{$U$}}
\put(56.00,7.00){\makebox(0,0)[cb]{$R_U$}}
\end{picture}
Let $\tau_R(U):=~\frac{|R_U|}{|U|}$ and $~\tau_L(U):=~\frac{|L_U|}{|U|}$.
The \emph{right thickness}  $\tau _{R}$ and \emph{left thickness} $\tau_{L}$ are
$$\tau_R(K):=~\inf_{U}\tau_R(U)~~~~~~~,~~~~~~~\tau_L(K):=~\inf_{U}\tau_L(U)$$
and the Newhouse thickness $\tau(K)$  is the minimum of lateral thicknesses $\tau_R(K)$ and $ \tau_L(K)$.

One of the special characteristics of above definition is the existence of an open and dense subset in
$C^{1+\epsilon}$ topology of regular Cantor sets whose elements have these lateral thicknesses varying continuously.
 But, Newhouse thickness is continuous in $C^{1+\epsilon}$ topology of all regular Cantor sets.

   Now we state some fundamental  results on stable intersection of regular Cantor sets:

I) If $HD(K)+HD(K')<1$, therefore no
translations of $K$ and $K'$ have
stable intersection (see [PT]),

II) If $\tau(K)\cdot\tau(K')>1$ and $K$ is linked to $K'$, then $(K,K')$ have stable intersection
(see [N] and [PT]),

III) If  $HD(K)+HD(K')>1$ and $K$ is linked to $K'$, then $(K,K')$  have generically stable intersection
(see [MY]).
\\Note that, $K$ is linked to $K'$ means closure of each gap of $K$ does not contain $K'$.

 Moreira  introduced affine Cantor sets with more than 2 expanding maps and small lateral thicknesses that having  stable intersection. He showed  that lateral thicknesses change continuously at affine Cantor sets defined by two expanding maps.
 Moreover, for affine Cantor sets defined by two affine expanding maps;

 $\bullet ~$ If $\tau_R(K)\cdot\tau_L(K')>1$, $\tau_L(K)\cdot\tau_R(K')>1$ and $K$ is linked to $K'$, then $(K,K')$ have stable intersection (see [M1]).

In [HMP], we constructed a pair of affine Cantor sets with the simplest possible combinatorics
which have stable intersection, while
  $\tau_R(K)\cdot\tau_L(K')<1$. These considerations motivate the following problem, that will be discussed in this work:
\begin{problem}
Does there exist a non empty open set in the space of affine Cantor sets defined by two expanding maps contained in the region
$$\big\{(K,K')~ |~~~\tau_R(K)\cdot\tau_L(K')<1~~ and~~ \tau_L(K)\cdot\tau_R(K')<1~\big\}$$
such that their elements have stable intersection?
\end{problem}

Indeed, we consider a pair of special  middle-$\alpha$ Cantor sets that  gives an affirmative solution to this problem. The main challenge is to construct a recurrent compact set of relative configurations, since, by the proposition in Subsection 2.3 of [MY], any relative configuration contained in a recurrent compact set is a configuration of stable intersection. The rest of this paper is outlined as follows:

In Section 2, we present the necessary definitions and state the recurrence condition on relative configurations of [MY] which implies stable intersection of pairs of regular Cantor sets.

In Section 3, we translate this condition to the setting of affine Cantor sets, which gives a recurrent condition on a simpler space of relative configurations.

Constructing a recurrent compact set of the relative configurations of
middle-$ \alpha$  Cantor sets, under special conditions, is the main subject of Section 4.

\section{ Basic definitions}

We will use notations similar to those of [MY] which are restated here.

Regular Cantor sets can also be defined as follows:

Let $A$ be a finite alphabet, $\mathcal{B}$ a subset of $A^{2}$, and
$\Sigma$ the subshift of finite type $A^{\mathbb{Z}}$ with allowed
transitions $\mathcal{B}$. We will always assume that $\Sigma$ is topologically mixing and every
letter in $A$ occurs in $\Sigma$.

An expanding map of type $\Sigma$ is a map $g$ with the following
properties:

{\rm i) the domain of $g$ is a disjoint union $\bigcup
_{\mathcal{B}}I(a,b)$, where for each $(a,b)$, $I(a,b)$ is a compact
subinterval of $I(a):=[0,~1]\times\{a\}$,

{\rm ii) for each $(a,b) \in \mathcal{B}$,
the restriction of $g$ to $I(a,b)$ is a smooth diffeomorphism
onto $I(b)$ satisfying $|Dg(t)|>1$ for all $t$.

The regular Cantor set associated to $g$ is the maximal invariant
set
\[
K:= \bigcap_{n\geq0}g^{-n}\big(\bigcup_{\mathcal{B}}I(a,b)\big).
\]
 These two definitions are equivalent. On one hand, we may, in the first definition, take $ I(i):=I_i$ for each $i \le r$, and, for each pair $i,j$ such that $\psi(I_i) \supset I_j$, take $I(i,j)=I_i \cap \psi^{-1}(I_j)$. Conversely, in the second definition, we can consider an abstract line containing all intervals $I(a)$ as subintervals, and $\{I(a,b)~|~ (a,b) \in \mathcal{B} \}$ as the Markov partition.

Also, a regular Cantor set $K$ is affine if $Dg$ be constant on every $I(a,b)$.

Let $\Sigma^{-}:=\big\{(\theta_{n})_{n\leq0}\mid~(\theta_{i},\theta_{i+1}%
)\in\mathcal{B}\ \emph{for}\ i<0~\big\}$. We equip $\Sigma^{-}$ with the following
ultrametric distance: for $\underline{\theta}\not =\underline{\tilde{\theta}%
}\in\Sigma^{-}$, set
\[
d(\underline{\theta},\underline{\tilde{\theta}}):=\left\lbrace
\begin{array}
[c]{lcr}%
1 &  & \theta_{0}\neq\tilde{\theta}_{0}\\
|I(\underline{\theta}\wedge\underline{\tilde{\theta}})| &  & \emph{otherwise}%
\end{array},
\right.
\]
where $\underline{\theta}\wedge\underline{\tilde{\theta}}:=(\theta_{-n}%
,\cdots,\theta_{0})$ if $\tilde{\theta}_{-j}=\theta_{-j}$ for $0\leq j\leq n$
and $\tilde{\theta}_{-n-1}\not =\theta_{-n-1}$.

Suppose that  $\underline{\theta}\in \Sigma^{-} $ and $n \in \mathbb{N}$, let $\underline{\theta}^n := (\theta_{-n}%
,\cdots,\theta_{0})$ and $ B(\underline{\theta}^n)$ be the affine map from $ I(\underline{\theta}^n)$ onto $ I({\theta}_0)$ which the diffeomorphism $ k \frac{\theta}{n}:= B(\underline{\theta}^n)\circ f _{\underline{\theta}^n}$ is orientation preserving. Then, for each $\underline{\theta}\in \Sigma^{-} $, there is a smooth diffeomorphism $ k^{\underline{\theta}}$ such that $ k \frac{\theta}{n}$ converges to $ k^{\underline{\theta}}$ in Diff$^{r}_{+}\big(I(\theta_{0}\big) $, for any $r \in (1, + \infty)$, uniformly in $ \underline{\theta}$.

Next, we define renormalization operators. For $(a,b)\in\mathcal{B}$, let
$$
f_{a,b}:=~[g|_{I(a,b)}]^{-1},
$$
this is a contracting diffeomorphism from $I(b)$ onto $I(a,b)$. If $\underline
a:=~(a_{0},a_{1},\ldots,a_{n})$ is a word of $\Sigma$, then we put
$$
f_{\underline a}:=~f_{a_{0},a_{1}} \circ\cdots f_{a_{n-1},a_{n}},
$$
this is a contracting diffeomorphism from $I(a_{n})$ onto a subinterval of
$I(a_{0})$, that we denote by $I(\underline a)$. Also let $F^{\underline\theta}$ be the affine map from $I(\theta_{0})$ onto $I(\theta_{-1},\theta_{0})$ with the same orientation of $f_{\theta_{-1},\theta_{0}}$.

Consider  $\mathcal{A}:=\big\{(\underline{\theta},A)\mid~\underline{\theta}%
\in\Sigma^{-}$ and $A:~I(\theta_{0})~\longrightarrow~\mathbb{R}$ is an affine embedding map\big\}. The renormalization operators $T_{\theta
_{1},\theta_{0}}:~\mathcal{A}~\longrightarrow~\mathcal{A}$ are defined by
\[T_{\theta_{1},\theta_{0}}(\underline{\theta},A):=~(\underline{\theta}\theta
_{1},A\circ F^{\underline{\theta}\theta_{1}}) ~~~~,~~~~ (\theta_{0},\theta_{1})\in\mathcal{B}.\]

For two sets of data $(A,\mathcal{B},\Sigma
,g)$ and  $(A^{\prime},\mathcal{B}^{\prime},\Sigma^{\prime},g^{\prime})$ defining
regular Cantor sets $K$ and $ K^{\prime}$, denote  $\mathcal{C}$ be the quotient of $\mathcal{A}\times\mathcal{A}$ by the diagonal action on the left of affine group.

A non empty compact set $\mathcal{L}$ in $\mathcal{C}$ is recurrent if for
every $u\in\mathcal{L}$ and $\ell,\ell^{\prime}\geq0$ with $\ell+\ell^{\prime}>0$,
when $(\underline{\theta},A)$ and  $(\underline{\theta^{\prime}},A^{\prime})$
represents $u$, then there exist words $\underline{a}=(a_{0},\cdots,a_{\ell}%
)$ and $\underline{a^{\prime}}=(a_{0}^{\prime},\cdots,a_{\ell^{\prime}}^{\prime})$
in $\Sigma$ and $\Sigma^{\prime}$, respectively, with $a_{0}=\theta_{0}$ and $a_{0}^{\prime}=\theta_{0}^{\prime}$,
such that $\big(T_{\underline{a}}(\underline{\theta},A),T_{\underline{a^{\prime}}}^{\prime}(\underline
{\theta^{\prime}},A^{\prime})\big)=v$ belongs to $\mathcal{L}^\circ:= \emph{int} \mathcal{L}$. The
following proposition has been proved in [MY].

\begin{proposition}
\noindent Any relative configuration (of limit geometries) contained in a
recurrent compact set is stably intersecting.
\end{proposition}

In the end of this section, we suppose that $\mathcal{S}:=\Sigma^{-}\times
\Sigma^{^{\prime}-}\times\mathbb{R^{\ast}}$, where $\mathbb{R^{\ast}}=\mathbb{R} \backslash \{0\}$. One
can see that the fibers of the quotient map $\mathcal{C}\longrightarrow
\mathcal{S}$ are one-dimensional and have a canonical affine structure. Moreover, this bundle map is trivializable.
We choose an explicit trivialization $\mathcal{C}\cong\mathcal{S}\times\mathbb{R}$ in order to have a coordinate in each fiber.

\section{ Transfer of renormalization operators on the space $\mathcal{S
}\times\mathbb{R}$}

Assume that $K$ is an affine Cantor set together with the Markov partition
$\{I(n,m)\}_{_{m,n\in A}}$ and an expanding  map
\[
\Phi\mid_{_{I(n,m)}}(x):=~p_{_{(n,m)}}\cdot x+q_{_{(n,m)}},
\]
where $A:=\{1,2,3,\cdots,N\}$. Let $I(n):=[a_{n}^{1}%
,~a_{n}^{2}]$ and $I(n,m):=[a_{n,m}^{1},~a_{n,m}^{2}]$, then for each
$(n,m)\in\mathcal{B}$, we have:
\[
\Phi\big(I(n,m)\big)=I(m).
\]
In the case $P_{(\theta_{0},\theta_{1})}>0$, map $F^{\underline{\theta}\theta_{1}}$ is
\[
F^{\underline{\theta}\theta_{1}}:I(\theta_{1})\longrightarrow I(\theta
_{0},\theta_{1}),
\]%
\[
F^{\underline{\theta}\theta_{1}}(x)=\frac{a_{\theta_{0},\theta_{1}}%
^{2}-a_{\theta_{0},\theta_{1}}^{1}}{a_{\theta_{1}}^{2}-a_{\theta_{1}}^{1}%
}(x-a_{\theta_{1}}^{1})+a_{\theta_{0},\theta_{1}}^{1}.
\]
Also in the opposite orientation, map $F^{\underline{\theta}\theta_{1}}$ is
\[
F^{\underline{\theta}\theta_{1}}(x)=\frac{a_{\theta_{0},\theta_{1}}%
^{2}-a_{\theta_{0},\theta_{1}}^{1}}{a_{\theta_{1}}^{2}-a_{\theta_{1}}^{1}%
}(x-a_{\theta_{1}}^{1})+a_{\theta_{0},\theta_{1}}^{2}.
\]
Therefore, in both cases, we obtain
\[
F^{\underline{\theta}\theta_{1}}(x)=\frac{1}{p_{(\theta_{0},\theta_{1})}%
}x-\frac{q_{(\theta_{0},\theta_{1})}}{p_{(\theta_{0},\theta_{1})}}.
\]

To continue, we  construct the homeomorphism between $\mathcal{S}%
\times\mathbb{R}$ and $\mathcal{C}$, then we transfer all renormalization
operators to $\mathcal{S}\times\mathbb{R}$.

\begin{theorem}
\noindent The map
\begin{align*}
L : ~\mathcal{C}~\longrightarrow~\mathcal{S}\times\mathbb{R}\\
\left[  (\underline{\theta},ax+b),(\underline{\theta}^{\prime},a^{\prime
}x+b^{\prime})\right]   & \mapsto(\underline{\theta},\underline{\theta
}^{\prime},\frac{a^{\prime}}{a},\frac{b^{\prime}-b}{a})
\end{align*}

is a homeomorphism between the space of relative configurations $\mathcal{C}
$ and $\mathcal{S}\times\mathbb{R}$.
\end{theorem}

\noindent\textbf{Proof.} $L$ is well defined:\newline Let
\[
\lbrack(\underline{\theta},ax+b),(\underline{\theta}^{\prime},a^{\prime
}x+b^{\prime})]=[(\underline{\theta},a_{1}x+b_{1}),(\underline{\theta}%
^{\prime},a_{1}^{\prime}x+b_{1}^{\prime})].
\]
Then there exist $c,d\in\mathbb{R}$ such that
\[
a_{1}x+b_{1}=c(ax+b)+d\ \ ,\ \ a_{1}^{\prime}x+b_{1}^{\prime}=c(a^{\prime
}x+b^{\prime})+d,
\]
therefore,
\[
\frac{a^{\prime}}{a}=\frac{a_{1}^{\prime}}{a_{1}}\ \ \ ,\ \ \frac{b^{\prime
}-b}{a}=\frac{b_{1}^{\prime}-b_{1}}{a_{1}}.
\]
$L$ is onto:
\[
\forall~(\underline{\theta},\underline{\theta}^{\prime},s,t)\in\mathcal{S}%
\times\mathbb{R},\hspace{1cm}L\Big([(\underline{\theta},x),(\underline{\theta
}^{\prime},sx+t)]\Big)=(\underline{\theta},\underline{\theta}^{\prime},s,t).
\]
Also, $L$ is one to one:\newline Suppose that
\[
L\Big([(\underline{\theta},x),(\underline{\theta}^{\prime},sx+t)]\Big)=L\Big([(\underline
{\theta},x),(\underline{\theta}^{\prime},s^{\prime}x+t^{\prime})]\Big),
\]
then $s=s^{\prime}$ and $t=t^{\prime}$.\newline From the structure of $L$, we see that
$L$ and $L^{-1}$ are continuous. $\Box$

Now we transfer the renormalization operators of relative configurations
$\mathcal{C}$ to the space $\mathcal{S}\times\mathbb{R}$. To do this, let
$\big((\theta_{1},\theta_{0}),(\theta_{1}^{\prime},\theta_{0}^{\prime}%
)\big)\in\mathcal{B}\times\mathcal{B}^{\prime}$, then we obtain
\begin{align*}
& (\underline{\theta},\underline{\theta}^{\prime},s,t)~ \xrightarrow{{L^{-1}}}
\ \left[  (\underline{\theta},x)~,~(\underline{\theta}^{\prime},sx+t)\right]  ~ \xrightarrow
{{ \big(T_{_{\theta_{1},\theta_{0}}},T_{_{\theta_{1}^{\prime},\theta_{0}^{\prime}}}^{\prime}\big) }} \\
& \left[  (\underline{\theta}~\theta_{1},\frac{1}{p_{(\theta_{0},\theta_{1})}%
}x-\frac{q_{(\theta_{0},\theta_{1})}}{p_{(\theta_{0},\theta_{1})}})\text{
},\text{ }(\underline{\theta}^{\prime}~\theta_{1}^{\prime},\frac{s}%
{p_{(\theta_{0}^{\prime},\theta_{1}^{\prime})}}x-\frac{q_{(\theta_{0}%
^{\prime},\theta_{1}^{\prime})}}{p_{(\theta_{0}^{\prime},\theta_{1}^{\prime}%
)}}s+t)\newline \right]  \overset{L}{\longrightarrow}\\
& \left(  \underline{\theta}~\theta_{1},\underline{\theta}^{\prime}%
~\theta_{1}^{\prime},\frac{p_{(\theta_{0},\theta_{1})}}{p_{(\theta_{0}%
^{\prime},\theta_{1}^{\prime})}^{\prime}}~s ~,~ p_{(\theta_{0},\theta_{1})}%
t-\frac{q_{(\theta_{0}^{\prime},\theta_{1}^{\prime})}^{\prime}}{p_{(\theta
_{0}^{\prime},\theta_{1}^{\prime})}^{\prime}} p_{(\theta_{0},\theta
_{1})}~s+q_{(\theta_{0},\theta_{1})}\right) .
\end{align*}

\bigskip We denote $L\circ\big(T_{_{\theta_{1},\theta_{0}}},T_{_{\theta
_{1}^{\prime},\theta_{0}^{\prime}}}^{\prime}\big)\circ L^{-1}$ by $T_{\big((\theta
_{1},\theta_{0}),(\theta_{1}^{\prime},\theta_{0}^{\prime})\big)} $ and for the sake of comfort we select those which are in these kinds:
\[
T_{\big((\theta_{1},\theta_{0}),id\big)}:~(\underline{\theta},\underline{\theta
}^{\prime},s,t)~\longrightarrow~(\underline{\theta}\theta_{1},~\underline
{\theta}^{\prime},~p_{_{\theta_{1},\theta_{0}}}s~,~p_{_{\theta_{1},\theta_{0}%
}}~t+q_{_{\theta_{1},\theta_{0}}}),
\]%
\[
T_{\big(id,(\theta_{1}^{\prime},\theta_{0}^{\prime})\big)}:~(\underline{\theta
},\underline{\theta}^{\prime},s,t)~\longrightarrow~(\underline{\theta
},~\underline{\theta}^{\prime}\theta_{1}^{\prime},~\frac{s}{p_{_{\theta
_{1}^{\prime},\theta_{0}^{\prime}}}^{\prime}}~,~t-\frac{q_{_{\theta_{1}%
,\theta_{0}}}^{\prime}}{p_{_{\theta_{1},\theta_{0}}}^{\prime}}s).~~~~~~
\]

\section{Construction of recurrent sets}
In this section, we introduce a pair of middle-$\alpha$ Cantor sets that have a recurrent compact
set in the  relative configurations, while  $\tau(K)\cdot
\tau(K^{\prime})<1$.
\begin{theorem}
Suppose that $K$ and $K^{\prime}$ are two homogenous Cantor sets with the convex hull $[0,~ 1]$ and expanding maps $\phi$  and  $\phi^{\prime}$ as

\unitlength1.30mm \linethickness{1pt} \begin{picture}(111.33,12.27)
\put(10.00,5.00){\makebox(0,0)[lc]{$K:$}}
\put(17.00,5.00){\line(1,0){9.00}}
\put(36.00,5.00){\line(1,0){9.00}}
\put(22.00,7.00){\makebox(0,0)[cb]{$\frac{1}{p}$}}
\put(14.00,7.00){\makebox(0,0)[cb]{$$}}
\put(54.00,5.00){\makebox(0,0)[lc]{$K':$}}
\put(41.00,7.00){\makebox(0,0)[cb]{$\frac{1}{p}$}}
\put(62.00,5.00){\line(1,0){10.00}}
\put(80.00,5.00){\line(1,0){10.00}}
\put(67.00,7.00){\makebox(0,0)[cb]{$\frac{1}{q}$}}
\put(77.00,7.00){\makebox(0,0)[cb]{$$}}
\put(85.00,7.00){\makebox(0,0)[cb]{$\frac{1}{q}$}}
\end{picture}
$$\phi(x):= \left\lbrace  \begin{array}
[c]{lcr}
~~~px ~~~~~~~~~~~~ x \in [0,~\frac{1}{p}]\\
px-p+1~~~~~  x \in [1-\frac{1}{p},~1]
\end{array} \right.
~~\phi^{\prime}(x):= \left\lbrace  \begin{array}
[c]{lcr}
~~~qx  ~~~~~~~~~~~~  x \in [0,~\frac{1}{q}]\\
qx-q+1~~~~~  x \in [1-\frac{1}{q},~1]
\end{array}, \right. $$
where
$$p:=\gamma^{40}:=(1.0321)^{40}=3.538923071...,$$
$$q:=\gamma^{31}:=(1.0321)^{31}=2.663024240....$$
Then pair $ (K , K')$ have stable intersection, while their  thickness product is smaller than one.
\end{theorem}

\noindent\textbf{Proof.} At first, we see that
$$\tau(K)\cdot\tau(K')=\frac{\frac{1}{p}}{1-\frac{2}{p}}\cdot\frac{\frac{1}{q}}{1-\frac{2}{q}}=
\frac{1}{(p-2)(q-2)}=0.980062299...$$

As  explained at the end of Section 3,  transferred renormalization operators can be considered
$$\mathbb{R^{\ast}\times \mathbb R~\longrightarrow~ \mathbb R^{\ast}\times \mathbb R}~~~~~~~~~~$$
$$(s,t)\overset{T_{0}}{\longmapsto}(ps~,~pt) ~~~~~~~~~~~~*~~~~~~~~~~~~(s,t)\overset{T_{1}}{\longmapsto}(ps~,~pt-p+1)$$
$$(s,t)\overset{T^{\prime}_{0}}{\longmapsto}(\frac{s}{q}~,~t)~~~~~~~~~~~~~~~~~~~~~~~~~~~~~
(s,t)\overset{T^{\prime}_{1}}{\longmapsto}(\frac{s}{q}~,~t+\frac{q-1}{q}s)$$

 Let $~s_{1}:=~\frac{\frac{1}{p}}{1-\frac{2}{q}+2(1-\frac{1}{q})\frac{1}{q^{39}}}$ ,
$~s_{2}:=~\frac{1-\frac{2}{p}}{\frac{1}{q}-2(1-\frac{1}{q})\frac{1}{q^{39}}}~$ and
$$ \Delta:=~\big\{(s,t)~|~~\gamma^{-1} s_{2} \leq s \leq \gamma s _{1}~~ ,~~-s+(1-\frac{1}{q})\frac {s}{q^{39}} \leq t\leq 1-(1-\frac{1}{q} )\frac{s}{q^{39}}~\big\}.$$
Note that, since
$~\tau_{1}\cdot\tau _{2} < 1 ~$and$~ (p-2)(q-2) < \gamma$, we have $ ~\gamma ^{-1}s _{2} <  s_{1} <s_{2} < \gamma s_{1} $.
\\Observe that the numerical approximations are

$~~~~~\gamma ^{-1}s _{2} \cong 1.1220161,~~s_{1} \cong 1.1349444,~~ s_{2}\cong1.1580329,~~ \gamma s_{1}\cong 1.1713761.~$

Take  $$\mathcal{L}:= ~\Delta\setminus\Delta_{1} \cup \Delta _{2}, $$ with
\\
$~~\Delta_{1}:=~\big\{(s,t)~ |~~
L^{1}:~t+(1- \frac{1}{q})(1+ \frac{1}{q^{39}})s > 1 ~~~~~ ~~ ~~,~~~
L^{2}: ~t+(1-\frac{1}{q} ) \frac{s}{q^{39}} > \frac{1}{p}
\\
~~~~~~~~~~~~~~~~~~~~~~~~~~~~~~~~~~~~,~~~L^3: ~t+(\frac{1}{q}-(1-\frac{1}{q})\frac{1}{q^{39}})s < 1-\frac {1}{p} ~\big\},$
\\$~~\Delta_{2}:= ~\big\{(s,t)~|~~
L^{4}:~ t+(1-(1-\frac{1}{q} ) \frac{1}{q^{39}}) s  < 1-\frac{1}{p}
 ~~~,~~~L^5:~ t+(\frac{1}{q}-(1-\frac{1}{q})\frac{1}{q^{39}})s< 0\\
~~~~~~~~~~~~~~~~~~~~~~~~~~~~~~~~~~~~,~~~L^6:~ t+(1- \frac{1}{q})(1+ \frac{1}{q^{39}})s >\frac{1}{p}~\big\}.
$

\begin{figure}[ht]
\centering 
\scalebox{0.5} 
{\includegraphics{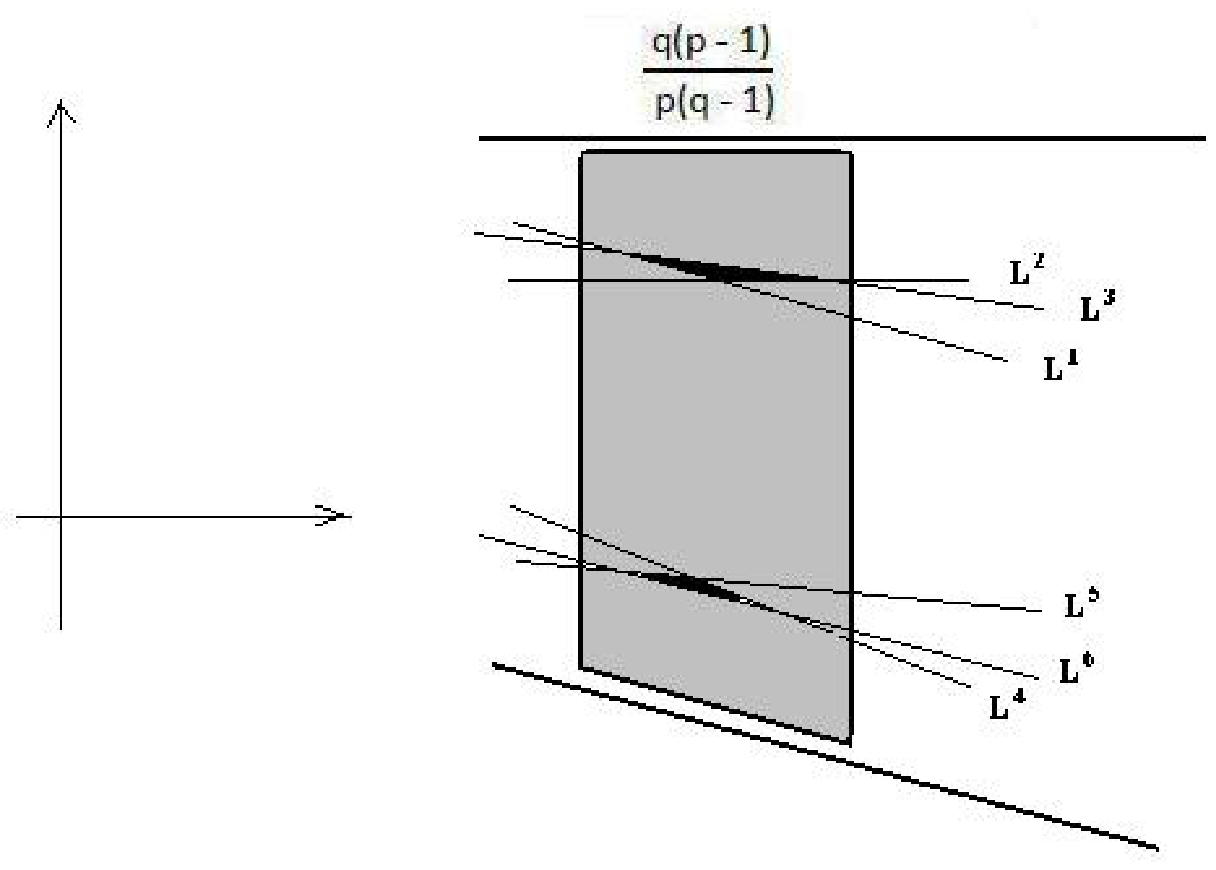}} 
\caption{}
\label{fig:exm} 
\end{figure}
We are going to show that  compact set $\mathcal{L}$, as displayed in Figure 1, is a recurrent set for the operators $(*)$.

Every vertical lines $s=s_{0}$~pass over itself with suitable compositions of the operators $(*)$ since $p^{31}=q^{40}$. Therefore  we can transfer the operators $(*)$ on these lines by
\begin{proposition}
 Let $s\in\mathbb R^{\ast}$  and  $\{a_{k}\}^{30}_{k=0}$ , $\{b_{k}\}^{39}_{k=0}$ be two finite sequences of numbers 0 and 1, then the maps
 \\ $~~~~~~~~~~T_s(t):=~p^{31}t+a_s ,      ~~~~~~~~~~~~~~~~~~~~~~~~~~~~~~~~~~~**~~~~~~~~~~~~  $
$~~~~~~~~~~~~~~~~~~~~~~~~~~~~~~~~~~~a_s:=-(p-1)p^{30}\big(\sum^{30}_{k=0}\frac{a_{k}}{p^{k}}-\frac{p(q-1)}{q(p-1)}s \sum^{39}_{k=0}\frac{b_{k}}{q^{k}}\big)~~~~$
\\are   return  maps to the vertical line $s$.
\end{proposition}

 \noindent\textbf{Proof}. Suppose that $\{b_{k}\}^{\infty}_{k=0}$  and $\{a_{k}\}^{\infty}_{k=0}$ are two arbitrary sequences of numbers 0 and 1. For every $a_k$ and $b_k$,  we can consider the operators  of $(*)$ in below form
 $$T _{a_{k}}(s,t):=~ (ps~ ,~ pt - (p-1) a_{k})~~~~,~~~~T _{b_{k}}(s,t):= ~(\frac{s}{q}~,~t+ (\frac{q-1}{q} )b_{k}s)$$
 Let $m,n\in \mathbb N$, then we obtain \\$~~~i)~ T _{a_{m-1}} \circ....\circ T_{a_{0}}(s,t)=\big(p^{m}s~,~ p^{m}t-(p-1)\sum^{m-1}_{k=0}a_{k}p^{m-1-k}\big),$
\\$~~~ii)~T _{b_{n-1}} \circ....\circ T_{b_{0}}(s,t)=\big(\frac{s}{q^{n}}~,~t+\frac{s}{q^{n}}(q-1)\sum^{n-1}_{k=0}b_{k}q^{n-1-k}\big).$

To prove the relations (i) and (ii), we use induction.  Case $ m=n=1$ is valid. Suppose that assertion satisfies for cases $i$ and $j$, then we have:
\\$~~~~T_{a_{i}}\circ T _{a_{i-1}} \circ....\circ T_{a_{0}}(s,t)=\big(p^{i+1}s~ ,~ p^{i+1}t-(p-1)\sum^{i-1}_{k=0}a_{k}p^{i-k}- (p-1)a_{i}\big)\\
~~~~~~~~~~~~~~~~~~~~~~~~~~~~~~~~~~~~~=\big(p^{i+1}s~,~ p^{i+1}t-(p-1)\sum^{i}_{k=0}a_{k}p^{i-k}\big),$\\
$~~~~T _{b_{j}}\circ T_{b_{j-1}} \circ....\circ T_{b_{0}}(s,t)=\big(\frac{s}{q^{j+1}}~,~t+ \frac{s}{q^{j}}(q-1)\sum^{j-1}_{k=0}b_{k}q^{j-1-k}+b_{j}\frac{(q-1)}{q}\cdot\frac{s}{q^{j}}\big)\\
 ~~~~~~~~~~~~~~~~~~~~~~~~~~~~~~~~~~~~~=\big(\frac{s}{q^{j+1}} ~,~t+\frac{s}{q^{j+1}}(q-1)\sum^{j}_{k=0}b_{k}q^{j-k}\big)$
\\and we see that the relations (i) and (ii) hold for cases $i+1$ and $j+1$.

 Replace  $m=31$ and  $n=40$ in the relations (i) and (ii), then we obtain
\\$~~~~~~~~~T _{b_{39}}\circ....\circ T_{b_{0}}\circ T _{a_{30}}\circ....\circ T _{a_{0}} (s,t)=\big(s~,~ p^{31}t- (p-1) \sum^{30}_{k=0}a_{k}p^{30-k}+\\ \frac{p^{31}s}{q^{40}}(q-1)\sum^{39}_{k=0}b_{k}q^{39-k}\big)=\big(s~,~ p^{31}t- (p-1)p^{30} ( \sum^{30}_{k=0}\frac{a_{k}}{p^{k}}-\frac{p(q-1)}{q(p-1)}s \sum^{39}_{k=0}\frac{b_{k}}{q^{k}})\big)$.

This completes the proof of proposition. $\Box$
\\If  $~\{a_{ik}\}^{30}_{k=0}~$  and  $~\{b_{ik}\}^{39}_{k=0}~$ be two arbitrary finite sequences of numbers $0, 1$ and   $s^{*}:=~\frac{p(q-1)}{q(p-1)}s$, then all of the  return   maps (or operators) are \\ $~~~~~~~~T_s(t)=p^{31}t+a_s ,      ~~~~~~~~~~~~~~~~~~~~~~~~~~~~~~~~~~~~~~~~~~~~~~~~~~~  ~~~~$
\\$~~~~~~~~~~~~~~~~a_s
:=-(p-1)p^{30}\big(\sum^{30}_{k=0}\frac{a_{ik}}{p^{k}}-s^{*} \sum^{39}_{k=0}\frac{b_{ik}}{q^{k}}\big)~
~~~~~~a_{ik},b_{ik}=0~or~1.$

Here we deal with $2^{31+40}$~squares of length $p^{-31}=\gamma^{-1240}=0.0{\underbrace{...}_{17}}0966...$, which project under angle $\theta:=\cot^{-1}s$~and determine this position, see Figure 2.
\begin{figure}[ht]
\centering 
\scalebox{0.3} 
{\includegraphics{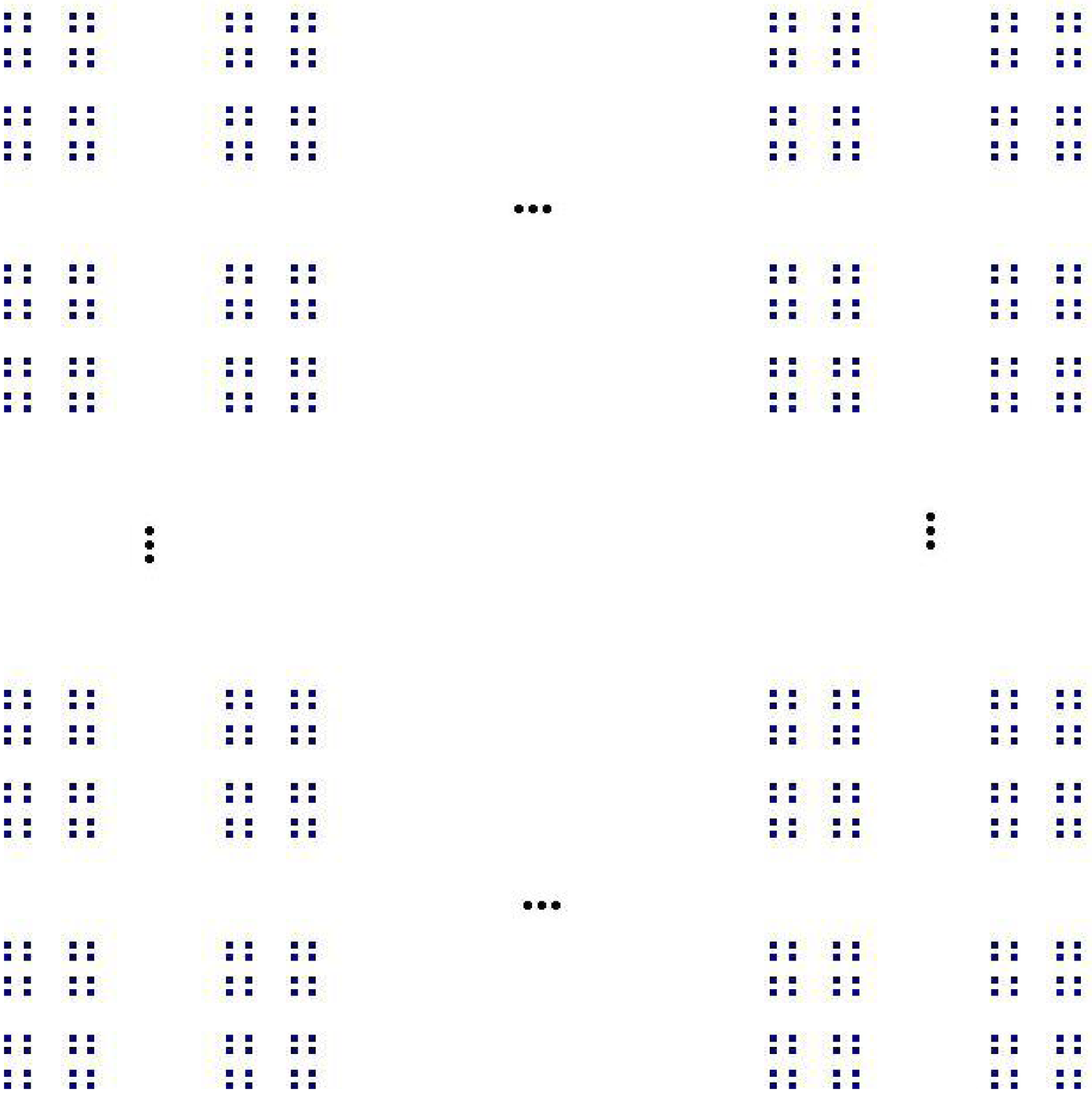}} 
\caption{}
\label{fig:exm} 
\end{figure}

\textbf{Definition.} The set $R\subset \mathbb{R} $ is a \emph{recurrent} if for every element of $R$ there exist suitable composites of maps $(**)$, so transfers that element to $R^{\circ}$.

 Let  $I_s^{\pm}:=~[a_{1},~b_{1}]\bigcup~[a^{\pm}_{2},~b^{\pm}_{2}]\bigcup~[a_{3},~b_{3}]$~ with
\\$~~~~~~a_{1}:=-\frac{q(p-1)}{p(q-1)}s^{*}+(1-\frac{1}{p})\frac{s^{*}}{q^{39}},~~~~~~~~~~~~~~~~~b_{1}:=
\frac{1}{p}-(1-\frac{1}{p})(1+ \frac{1}{q^{39}})s^{*},$\\
$~~~~~~a_{2}^{-}:=-\frac{p-1}{p(q-1)}s^{*}+(1-\frac{1}{p})\frac{s^{*}}{q^{39}},~~~~~~~~~~~~~~~~b_{2}^{-}:=
1-(1-\frac{1}{p})(1+\frac{1}{q^{39}})s^{*},$
\\$~~~~~~a_{2}^{+}:=1-\frac{1}{p}-\frac{q(p-1)}{p(q-1)}s^{*}+(1-\frac{1}{p})\frac{s^{*}}{q^{39}},~~~~
~~~~b_{2}^{+}:=\frac{1}{p}-(1-\frac{1}{p})\frac{s^{*}}{q^{39}},$
\\$~~~~~~a_{3}:=1-\frac{1}{p}-\frac{p-1}{p(q-1)}s^{*}+(1-\frac{1}{p})\frac{s^{*}}{q^{39}},~~~~~~~~~b_{3}:=
1-(1-\frac{1}{p})\frac{s^*}{q^{39}}.$
\\
\\\put(10.00,1.00){\makebox(0,0)[lc]{$I_s^{\pm}:~$}}
\put(20.00,1.00){\line(1,0){20.00}}
\put(42.00,1.00){\line(1,0){20.00}}
\put(64.00,1.00){\line(1,0){20.00}}
\put(22.00,2.00){\makebox(0,0)[cb]{$a_1$}}
\put(39.00,2.00){\makebox(0,0)[cb]{$b_1$}}
\put(44.00,2.00){\makebox(0,0)[cb]{$a_2^{\pm}$}}
\put(61.00,2.00){\makebox(0,0)[cb]{$b_2^\pm$}}
\put(65.00,2.00){\makebox(0,0)[cb]{$a_3$}}
\put(84.00,2.00){\makebox(0,0)[cb]{$b_3$}}
\\
\begin{lemma} For every  $s_1\leq s\leq \frac{q(p-1)}{p(q-1)}~$, the set $I_s^{-}$ and for every  $\frac{q(p-1)}{p(q-1)}<s\leq s_2,~$ the set  $I_s^{+}$  is a  recurrent set for maps $(**)$.
\end{lemma}
\textbf{Proof}.  If $s_1\leq s \leq s_2$, then
\\$~~~~~~~\delta _{1} := 0.987915116...=\frac{q-1}{(p-1)(q-2)+ \frac{2(p-1)(q-1)}{q^{39}}}<s^{*}\\~~~~~~~~~~~
~~~~~~~~~~~~~~~~~~~~~~~~~~~~< \frac{(p-2)(q-1)}{(p-1)-\frac{2(p-1)(q-1)}{q^{39}}} = 1.00801257... =:\delta_{2}. $

 Fix $s= cotg~\theta$ and relinquish of  notation $\pm$  on our calculations.

 At first, we show that for every point in interval
$$J_s:=~\big[0-\frac{q(p-1)}{p(q-1)}(1-\frac{1}{q}+\frac{1}{q^{2}})\frac{s^{*}}{q^{38}}~,
~~\frac{1}{p^{29}}-\frac{q(p-1)}{p(q-1)}(1-\frac{1}{q})\frac{s^{*}}{q^{39}}\big],$$
there exist suitable maps of $(**)$ that send that point to $I_s^{\circ}$.
We  remind that $J_s$ is projection of the union of these 14 squares

$$ C_{ij}:=~[0,~\frac{1}{p^{31}}]\times[0,~\frac{1}{q^{40}}]+ {(c_{i},c_{j})}~~, ~~1\leq i,j \leq 4$$
with
\\$(c_{i},c_{j})\in\{c_{i}\}^{4}_{i=1}\times\{c_{j}\}^{4}_{j=1}:=~\big\{0~,~(1-\frac{1}{p})\frac{1}{p^{30}}~,~
(1-\frac{1}{p})\frac{1}{p^{29}}~,~(1-\frac{1}{p})(\frac{1}{p^{29}}+\frac{1}{p^{30}})\big\}\times$\\ $~~~~~~~~~~~~~~~~~~~~~~~~~~~~~~~~~~~~~~~~~\big\{0~,~(1-\frac{1}{q})\frac{1}{q^{39}}~,~(1-\frac{1}{q})\frac{1}{q^{38}}~,~
(1-\frac{1}{q})(\frac{1}{q^{38}}+\frac{1}{q^{39}})\big\},$
\\except the squares~~$ C_{41}$ and $ C_{14}.$ Let C be one of these squares, $ T $ be its correspondent operator and $\Pi_{\theta}$ be the projection map on $\mathbb{R}^{2}$. Set $K(C):=~\Pi_{\theta}(C) -T^{-1}(I_s^{\circ})$ contains  4 components if $s\in (s_1,~s_2) $ and 2 components if $s=s_1$ or $s=s_2$.

We show that intervals $\Pi_{\theta}(C)$ overlap each other and each point of $K(C)$ goes to $I_s^{\circ}$ under other operators of $(**)$. Suppose that $T_{ij}(t)=p^{31}t+a_{ij}$ and $T_{i'j'}(t)=p^{31}t+a_{i'j'}$ are the corresponding  operators to the squares $C_{ij}$ and $C_{i'j'}$, respectively.
Therefore we have $ T_{i'j'}(t)=T_{ij}(t)-(a_{ij}-a_{i'j'})$ and we see that on special conditions point t goes to
$I^{\circ}$ under the operator $ T_{i'j'}$, in fact,
\\ $~~~~i)$ If \\ $~~~~~max~~\big\{-\frac{1}{p}+\frac{(p-1)(q-2)}{p(q-1)}s^{*}+2(1-\frac{1}{p})\frac{s^{*}}{q^{39}}~~,~~
 1-\frac{2}{p}-\frac{(p-1)}{p(q-1)}s^*+2(1-\frac{1}{p})\frac{s^{*}}{q^{39}}\big\}\\~~~~~~~~~~~~= max~~a_{2}^{\pm}-b_{1}<a_{ij}-a_{i'j'}
 <\frac{1}{p}+\frac{p-1}{p(q-1)}s^{*}-2(1-\frac{1}{p})\frac{s^{*}}{q^{39}}=b_{1}-a_{1},~~~$
\\then all points  of  last 3 components of $K(C_{ij})$ will be sent to $I_s^{\circ}$ by the operator $T_{i'j'}$ and right side of $\Pi_{\theta}(C_{ij})$ intersects left side of $\Pi_{\theta}(C_{i'j'}),$
\\ $~~~~ ii)$  If \\$~~~~~~~~~~~1-b_{1}=1-\frac{1}{p}+(1-\frac{1}{p})(1+\frac{1}{q^{39}})s^{*}<a_{ij}-a_{i'j'}\\
~~~~~~~~~~~~~~~~~~~~~~~~~~~~~~~~~~~~~~~~~~~~~~<1+\frac{q(p-1)}{p(q-1)}s^{*}-2(1-\frac{1}{p})\frac{s^{*}}{q^{39}}
=b_{3}-a_{1},$
\\then the last component of $K(C_{ij})$ will be sent to $I_s^{\circ}$ by the operator $T_{i'j'} $ and right side of  $ \Pi_{\theta}(C_{ij})$ intersects left side of $\Pi_{\theta}(C_{i'j'}). $

Similar result is valid for negative cases.

We need numeral values of right and left sides of above inequality, indeed,
\\$~~~max~~ a_{2}^{\pm}-b_{1}= 0.008670035...,$
$~~~~~~~min ~~b_{1}-a_{1}= 0.708758125...,$
\\$~~~max~~~~ 1-b_{1}= 1.440604766...,$
$ ~~~~~~~min~~ b_{3}- a_{1}= 2.134944413....$

We split interval $J_s$ in such cases:
\\ \textbf{Case 1.}  $~t\in K (C_{31}) $, we consider the operators $T_{41}$ and $T_{42}$  as
\\$ ~~~~~T_{41}(t)= p^{31}t-(p-1)p^{30}(-\frac{s^{*}}{q^{38}}-\frac{s^{*}}{q{39}})=p^{31}t+ \frac{p-1}{p}(q^{2}+q)s^{*},$
\\$ ~~~~~T_{42}(t)= p^{31}t-(p-1)p^{30}(\frac{1}{p^{30}}-\frac{s^{*}}{q^{38}}-\frac{s^{*}}{q{39}})=p^{31}t+\frac{p-1}{p}(-p+(q^{2}+q)s^{*}).$
\\ On the other hand, the corresponding operator    $C_{31} $ is
\\$~~~~~T_{31}(t)= p^{31} t -(p-1)p^{30}(-\frac{s^{*}}{q^{38}})=p^{31}t+\frac{p-1}{p}(q^{2})s^{*}$
\\and we obtain
\\ $~~~~~~~~~~~~~~~-1.925836829< a _{31}-a_{41}=\frac{p-1}{p}(-q)s^{*}<-1.887440068,$
\\$ ~~~~~~~~~~~~~~~+0.613086242 < a_{31}-a_{42} = \frac{p-1}{p}(p-qs^{*})<+0.651483003,$
\\ therefore the element t will be  sent to $I_s^{\circ}$ by using the operators $T_{41}$ or $T_{42}$, and
 interval $\Pi_{\theta}(C_{31})$ connects $\Pi_{\theta}(C_{41})$ to $\Pi_{\theta}(C_{42}).$
\\ \textbf{Case 2.} $~t\in K(C_{42})$, here we need the operator
\\ $~~~~~T_{32}(t)=p^{31}t-(p-1)p^{30}(\frac{1}{p^{30}}-\frac{s^{*}}{q^{38}})=p^{31}t+\frac{p-1}{p}(-p+q^{2}s^{*})$
\\and we obtain \\ $~~~~~~~~~~~~~~~~~~~~~~~~~a_{42}-a_{32}=\frac{p-1}{p}(q)s^{*}=a_{41}-a_{31},$\\
therefore  first 3 components of $K(C_{42})$ map to $I_s^{\circ}$ by the operator $T_{31}$
 since \\$~~~~~~~~~~~~~~~~~~~~~~- 0.708758125...<a_{42}-a_{31}<-0.008670035...$
\\and the last component of $K(C_{42})$ maps to $I_s^{\circ} $ by using the operator $T_{32}$ and interval $\Pi_{\theta}(C_{42})$ connects $\Pi_{\theta}(C_{31})$ to $\Pi_{\theta}(C_{32}).$
\\ \textbf{Case 3.} $~t\in K(C_{32})$, in this case, we consider the operator
\\ $~~~~~T_{21} (t) = p^{31}t-(p-1)p^{30}(-\frac{s^{*}}{q^{39}})=p^{31}t+\frac{p-1}{p}(q)s^{*} $
\\ and we have:\\ $ ~~~~~~~~~~+0.599935514<a_{32}-a_{21}=\frac{p-1}{p}(-p+(q^{2}-q)s^{*})<+0.663790257,$
\\ therefore the operators $T_{42}$ or $ T_{21}$ send the point t to $I_s^{\circ}$ and interval $\Pi _{\theta}(C_{32})$ connects $\Pi(C_{42})$ to $\Pi_{\theta}(C_{21}).$
\\ \textbf{Case 4.} $~t\in K (C_{21})\bigcup K (C_{22}) \bigcup K (C_{11})$, when $ t\in K (C_{21})$, we have  $T_{11} (t)= p^{31}t $ and we see that numbers $a_{21}-a_{32}$ and $ a_{21}-a_{11} $ satisfy  $(i)$ and $(ii)$ conditions above, this says that  point t goes to $I_s^{\circ}$ under the operators $T_{32}$ or $T_{11}$ and interval $\Pi_{\theta} (C_{21})$ connects $\Pi_{\theta}(C_{32})$ to $\Pi(C_{11}).$ Elements $K(C_{22}) $ and $K(C_{11})$ behave like  (2) and (3) cases above.
\\ \textbf{Case 5.} $~t\in K(C_{12} )\bigcup K (C_{43})$, in this case, we obtain\\
$~~~~~~~~~~~~~~~a_{22} = \frac{p-1}{p} (-p+qs^{*}), ~~~~~~~~~~~~~~~~~a_{12}=\frac{p-1}{p}(-p)$,\\
$ ~~~~~~~~~~~~~~~a_{43}=\frac{p-1}{p}(-p^{2}+(q^{2}+q)s^{*}),~~~~~~~a_{33}=\frac{p-1}{p}(-p^{2}+q^{2}s^{*}),$
\\therefore,
\\ $~~~~~~~~~~~~~~~~~~~~~~~~~~~~a_{12}-a_{22}=\frac{p-1}{p}(-q)s^{*}=a_{33}-a_{43},$\\
   $~~~-0.608256623<a_{12} -a_{43}=\frac{p-1}{p} (p^{2}-p-(q^{2}+q)s^{*})<-0.467608362,$ \\as above discussion, the  point t goes to set $I_s^{\circ}$ and the projections cover together.
\\ \textbf{Case 6.} $~$Other cases, here we use  relation $$T_{ij}(t)=T_{i(j-2)} \big(t-(\frac{p-1}{p})\frac{1}{p^{29}}\big)$$ that sends the point $t$ to $I_s^{\circ} $  since point $ t-(\frac{p-1}{p})\frac{1}{p^{29}} $ comes back to above cases.

Till now, we have shown that points  of interval $J_s$ will be sent to $I_s^{\circ}.$

Let $D$ be the union of these 14 squares and  $D_{1} $ be all the squares in  $ [ 0, ~\frac{1}{p^{28}}]\times [0,~ \frac{1}{q^{37}}] $, except two corner squares $[ 0,~ \frac{1}{p^{31}}]\times[\frac{1}{q^{37}}-\frac{1}{q^{40}},~\frac{1}{q^{37}}]$ and $ [\frac{1}{p^{28}}-\frac{1}{p^{31}},~ \frac{1}{p^{28}} ] \times [0, ~\frac {1}{q^{40}}]$. We are going to pave all the squares in $[0,~\frac{1}{p}]\times[0,~\frac{1}{q}]$ by set $D_{1} $, of course, except two corner squares $$[0,~\frac{1}{p^{31}}]\times[\frac{1}{q}-\frac{1}{q^{40}},~\frac{1}{q}]~~~ ,~ ~~ [\frac{1}{p}-\frac{1}{p^{31}},~\frac{1}{p}]\times[0,~ \frac{1}{q^{40}}].$$
We represent this set to $\mathfrak{D} $ and we construct it by induction.

Let $D_{n}$ be the set of all the squares in $[0,~\frac{1}{p^{i+1}}]\times[0,~\frac{1}{q^{j+1}}]$, except two corner squares\\
$~~~~~~~~~~~~~~~~[0,~\frac{1}{p^{31}}]\times[\frac{1}{q^{j+1}}-\frac{1}{q^{40}},~\frac{1}{q^{j+1}}]~~~,~~~ [\frac{1}{p^{i+1}}-\frac{1}{p^{31}},~ \frac{1}{p^{i+1}}]\times[0,~\frac{1}{q^{40}}],$\\ for a pair $(i,j)$ with condition $~1\leq i \leq 28$,~ $1\leq j\leq 37$ and also $D_{n}$ supports that every line with slope  $\theta=\theta_s$, that passes among  itself, has a conflict with a copy of $D_1$.  Before  defining  set $D_{n+1} $, as shown in Figure 3, put
 $$E_{n}: =~D_{n}+\big((1-\frac{1}{p})\frac{1}{p^{i}} , {0}\big)~~~~~,~~~~~F_{n}:=~D_{n}+\big(0,(1-\frac{1}{q})\frac{1}{q^{j}}\big)$$

\begin{figure}[ht]
\centering 
\scalebox{0.3} 
{\includegraphics{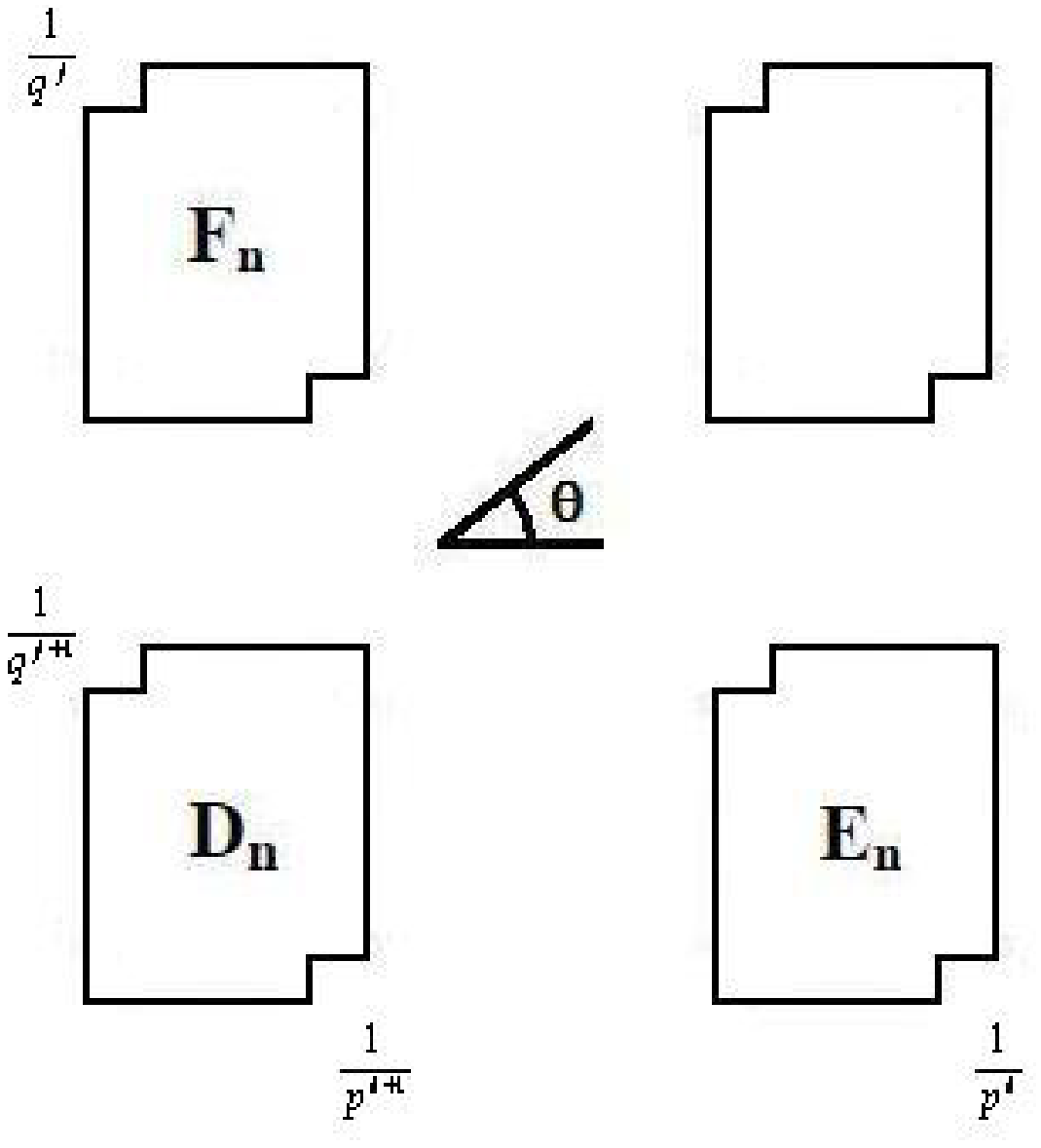}} 
\caption{}
\label{fig:exm} 
\end{figure}

I) If $\Pi_{\theta}(D_{n}) \bigcup \Pi_{\theta} (E_{n})$ be connected, then we let
$ D_{n+1} $ be the set of all the squares in \\$[0,~\frac{1}{p^{i}}]\times[0,~\frac{1}{q^{j+1}}]$, except $[0,~\frac{1}{p^{31}}]\times[\frac{1}{q^{j+1}}-\frac{1}{q^{40}},~\frac{1}{q^{j+1}}]$ and $[\frac{1}{p^{i}}-\frac {1}{p^{31}},~ \frac{1}{p^{i}}]\times[0,~\frac{1}{q^{40}}]$

otherwise,

II) If $\Pi_{\theta}(D_{n}) \bigcup \Pi_{\theta}(F_{n}) $ be connected, then we let
 $ D_{n+1}$ be the set of all the squares in \\$[0,~\frac{1}{p^{i+1}}]\times[0,~\frac{1}{q^{j}}]$, except
$[0,~\frac{1}{p^{31}}]\times[\frac{1}{q^{j}}-\frac{1}{q^{40}},~\frac{1}{q^{j}}]$ and $[\frac{1}{p^{i+1}}-\frac {1}{p^{31}},~ \frac{1}{p^{i+1}}]\times[0,~\frac{1}{q^{40}}]$.

Condition (I) is equivalent with  the relation
\\ $ (1-\frac{1}{p})\frac{1}{p^{i}} - \Big(\frac{1}{q^{j+1}}-(1-\frac{1}{q})\frac{1}{q^{39}}\Big)\frac{q(p-1)}{p(q-1)}s^{*}<\frac{1}{p^{i+1}} -(1-\frac{1}{q} )\frac{1}{q^{39}}\cdot\frac {q(p-1)}{p(q-1)}s^{*}\\ \xRightarrow {\times \frac{p^{i+1}}{s^{*}}}
\frac{p-2}{s^{*}} < \frac{p^{i}}{q^{j}} (\frac{p-1}{q-1}) -\frac{2(p-1)p^{i}}{q^{39}} \Longrightarrow \frac {(p-2)(q-1)}{(p-1)s^{*}}+\frac{2(q-1)p^{i}}{q^{39}}<\frac{p^{i}}{q^{j}}.~~~~~~~~~~~~~~~~(1) $

Condition (II) is equivalent with the relation
\\$ 0-\Big(\frac{1}{q^{j+1}}-(1-\frac{1}{q})\frac{1}{q^{39}}\Big)\frac{q(p-1)}{p(q-1)}s^{*}< \frac{1}{p^{i+1}} -(1-\frac{1}{q}) (\frac{1}{q^{j}}+ \frac{1}{q^{39}})\frac {q(p-1)}{p(q-1)}s^{*}\\ \xRightarrow {\times \frac{q-1}{(p-1)s^{*}}p^{i+1}} \frac {p^{i}}{q^{j}}< \frac{q-1}{(p-1)(q-2)s^{*}}- \frac {2(q-1)p^{i}}{(q-2)q^{39}}.~~~~~~~~~~~~~~~~~~~~~~~~~~~~~~~~~~~~~~~~~~~~~(2)$

At first, we show that assertion is valid for $n=1$. We can put $i=28$ and $j=37$ in the relation (2) and then we obtain that interval $\Pi_{\theta}\big(D+(0,(1-\frac{1}{q})\frac{1}{q^{37}})\big)$ intersects interval $ \Pi_{\theta}(D)$. This matter is similar for sets
$\Pi_{\theta}\big(D+((1-\frac{1}{p})\frac{1}{p^{28}},0)\big)$ and
$\Pi_{\theta}\big(D+((1-\frac{1}{p})\frac{1}{p^{28}},(1-\frac{1}{q})\frac{1}{q^{37}})\big)$.
The relation (1) is not valid here, and we had to use their operators. Two squares that occur on points
$$ (c_{1},c_{4})=\big((1-\frac{1}{p})(\frac{1}{p^{29}}+ \frac {1}{p^{30}}), 0\big)
~, ~(c_{8},c_{5}):=~ \big((1-\frac{1}{p})\frac{1}{p^{28}} , (1-\frac{1}{q})(\frac{1}{q^{37}}+\frac{1}{q^{38}}+\frac{1}{q^{39}})\big)$$
have the following operators
\\$~~~~~T_{14}(t) =p^{31}t+\frac{p-1}{p}(-p^2-p),$
\\$~ ~~~~T_{85}(t) = p^{31}t+\frac{p-1}{p} \big(-p^{3}+(q^{3}+q^{2}+q)s^{*}\big)$
\\and we get
\\$~~~~~+0.025457501<a_{85}-a_{14} =\frac{p-1}{p}\big(-p^{3}+p^{2}+p+(q^{3}+q^{2}+q)s^{*}\big)<0.438403979.$

Therefore the projection of  two squares $C_{84}$ and $C_{14}$ fill each other's  middle gaps, also the square $C_{84}$ set upper than
 $C_{14}$ under  projection $\Pi_{\theta}$, this expresses that  boundary points fall in set $J_s$ or
$J_s+\Pi_{\theta}\big((1-\frac{1}{p})\frac{1}{p^{28}},(1-\frac{1}{q})\frac{1}{q^{37}})\big).$

Suppose that assertion is valid for $n$ and  we check property of  $D_{n+1}$.
\\Since $(p-2)(q-2) > 1 $,  the relations (1) and (2) do not happen together, therefore  the process stops when
 $$ \frac{q-1}{(p-1)(q-2)s^{*}} - \frac{2(q-1)p^{i}}{(q-2)q^{39}}-1 < \frac{p^{i}}{q^{j}}-1< \frac {(p-2)(q-1)}{(p-1)s^{*}} +\frac{2(q-1)p^{i}}{q^{39}} -1 .$$

In fact, we show that always
\\$ ~~~~~~~|\frac{p^{i}}{q^{j}} - 1 |>max \Big\{\frac {(p-2)(q-1)}{(p-1)\delta_1}+\frac{2(q-1)p^{i}}{q^{39}}-1~,~-\big(\frac {q-1}{(p-1)(q-2)\delta_2} - \frac{2(q-1)p^{i}}{(q-2)q^{39}}-1\big)\Big\}
\\ ~~~~~~~=max \big\{0.020343299...+ \frac{2(q-1)p^{i}}{q^{39}}~,~0.019937701... +\frac{2(q-1)p^{i}}{(q-2)q^{39}}\big\}.~~~~~~~~~~(3)$

To prove the relation (3), we divide $1\leq i \leq 27$ in the following cases:
\\   1) If $i=27 $, then $\frac{p^{i}}{q^{39}} =\frac{q}{p^{4}}$ and we obtain
\\$~~~~~~~ min \big\{|\frac{p^{27}}{q^{j}}-1|~\mid~~1\leq j \leq 37 ~\big\} =~min \big\{|\frac{p^{27}}{q^{34}}-1|~,~ |\frac{p^{27}}{q^{35}}-1| \big\}\\   ~~~~~~~=~ min\{|\gamma^{26} - 1|~,~|\gamma^{-5} -1|\} = 0.146131269.... $

On the other hand,
 \\$~~~~~~\frac {2(q-1)p^{27}}{q^{39}}=\frac {2(q-1)q}{p^{4}} = 0.056470184... ~~~,~~~\frac {2(q-1)p^{27}}{(q-2)q^{39}}=\frac {2(q-1)q}{(q-2)p^{4}} =0.085170618...,  $
 \\ therefore the relation (3) happens anytime.
\\ 2)  If $i=26$, then $ \frac {p^{26}}{q{39}} = \frac{q}{p^{5}}$ and by
\\$ ~~~~~~~~min \big\{|\frac {p^{26}}{q^{j}} -1|~|~~ 1\leq j\leq 37 ~\big\}= min \big\{|\frac{p^{26}}{q^{33}}-1|~,~ |\frac{p^{26}}{q^{34}}-1| \big\} \\ ~~~~~~=~min\{|\gamma^{17} -1|~,~|\gamma^{-14} - 1|\} = 0.357467488...$
\\the  relation(3) is valid.
\\ 3) If $ i\leq 25$, then maximum value in the right side of the  relation (3) happens when $i=25$, thus,
$$  \frac{2(q-1)p^{i}}{q^{39}} = 0.004508966...~~~,~~~\frac{2(q-1)p^{i}}{(q-2)q^{39}}=0.006800605....$$

We show that $|\frac {p^{i}}{q^{j}}-1| > 0.028738306$, where $1\leq j\leq 37$.
Since $\frac{p^{i}}{q^{j}}=\gamma ^{40i-31j}$ and $40i -31j \neq 0$, value $$min \big\{| \frac{p^{i}}{q^{j}}-1|~|~~ 1\leq j\leq37 ~~,~~ 1\leq i \leq25 ~~\big\}$$
 happens when $40i-31j=\pm 1$, (for example $i=7$ and $j=9$).
Also, we have:
 $$min \{|\gamma-1| ~,~|\gamma ^{-1} -1|\}=0.031101637....$$

Therefore  $ D_{n+1}$  always exists since sets $E_{n}$ and $F_{n}$ are copies of $D_{n}$ and the relation (3) is valid. This  process continues until we construct  $\mathfrak{D}$, here is $i=1 $ and $j=1$ since $ (p-2)(q-2)> 1 $. Therefore, every line $L_{\theta}$ that passes among the convex hull $\mathfrak{D}$, had to meet at least one copy of set $D_{1}$  and always projected point will be sent to $I_s^{\circ}$ by using one of the operators $(**)$, like case (6) above. Also,  sets $\mathfrak{D}+(1- \frac{1}{p}, 0),$  $\mathfrak{D}+(0,1- \frac{1}{q} )$ and
 $\mathfrak{D}+(1- \frac{1}{p},1- \frac{1}{q} )$ behave like set $\mathfrak{D}$, as  case (6) above.

For $s=s_1$,  $I_s$ is  connected interval $(-s,~1)$,  but  our images are not the points which have been obtained from the intersection of lines $L^3 , L^1$ and $L^6 , L^5$ yet. The same  situation holds for $ s=s_{2}.$

Here the  proof is finished. $\Box$\\

\textbf{Lemma 2} . If $ \gamma ^{-1} s_{2}\leq s \leq s_{1}$ or $s_{2}\leq s \leq \gamma   s_{1}$, then interval $I_{s} =(-s,~ 1) $ is a recurrent set for the maps $(**)$.

\textbf{Proof}. It is enough to show that every line with slope $ \theta= \cot ^{-1} s $ conflicts with one of the squares in $[0,~1] \times [0,~1]$. Like the above lemma, we do this by induction.

Let $D_{1}:=C_{11}$ and suppose that $ D_{n}$ be the set of all the squares in $[0,~\frac{1}{p^{i+1}}] \times [0,~\frac{1}{q^{j+1}}]$, for a pair of $i,j $ and every line $L_{\theta }$ intersects one of the squares $D_{n}$. Let $E_{n}, F_{n}$ be like Lemma 1 and

I) Let $ D_{n+1} $ be the set of all the squares in $[0,~\frac{1}{p^{i}}]\times[0,~\frac{1}{q^{j+1}}]$ if $\Pi_{\theta}(D_{n}) \bigcup \Pi_{\theta} (E_{n})$ be connected otherwise,

II) Let  $ D_{n+1}$ be the set of all the squares in $[0,~\frac{1}{p^{i+1}}]\times[0,~\frac{1}{q^{j}}]$ if $\Pi_{\theta}(D_{n}) \bigcup \Pi_{\theta}(F_{n}) $ be connected.
\\If we take  small sentences in the relations (1) and (2) of previous lemma  and replace  $s^{*}:=~\frac{p(q-1)}{q(p-1)}s$, then the process stops when
\\$ ~~~~~~~~~~~~~~~~~~~~~~~~~~~~~~~~~~~~~\frac {q}{ p(q-2)s} <\frac {p^{i}} {q^{j}} <\frac {q}{ps}(p-2).$
\\We need to show that the relation $ s_{1} < \frac {p^{i}}{q^{j}}  s < s_{2}$ does not happen any time:
\\1)~If  $\gamma ^{-1} s_{2}\leq s \leq s_{1}$, then\\
$~~~~~~~~~~~~~~~~~~~~p^{i}> q^{j}~ \Longrightarrow ~1< \gamma < \frac {p^{i} } {q^{j}} ~\xRightarrow {\gamma ^{-1} s_{2} <s}~ s_{2} <\frac{p^{i}}{q^{j}} s, $
\\$ ~~~~~~~~~~~~~~~~~~~~p^{i}< q^{j} ~ \Longrightarrow ~ \frac{p^{i}}{q^{j}}< \gamma^{-1}  <1 ~\xRightarrow {s<s_{1}} \frac{p^{i}}{q^{j}}~ s < \gamma ^{-1} s_{1} < s_{1}.$
\\2)~If $ s_{2} \leq s\leq \gamma s_{1} $, then
\\$~~~~~~~~~~~~~~~~~~~~~ p^{i} > q^{j}~ \Longrightarrow ~1< \gamma < \frac {p^{i}}{q^{j}}~\xRightarrow {s_{2} < s }~  s_{2}< \gamma s _{2} < \frac{p^{i}} {q^{j}} s,$
\\$~~~~~~~~~~~~~~~~~~~~~p^{i}< q^{j} ~ \Longrightarrow ~ \frac{p^{i}}{q^{j}}< \gamma^{-1}  <1 ~\xRightarrow {s< \gamma s_{1}}~\frac{p^{i}}{q^{j}} s <  s_{1}.$
\\Thus,  construction continues till $i=j=0$. This completes the proof of the lemma. $\Box$

Now we show that $\mathcal{L}$ is a recurrent set for the operators $(*)$.

\textbf{Case1.} If   $~s_1\leq s \leq s_2$, then we use  Lemma 1.

\textbf{Case2.} If $~\gamma ^{-1} s_{2}< s < s_{1}$ or $s_{2}< s < \gamma   s_{1}$, then  we use  Lemma 2 and the  operators \\$~~~~~~~~\pi_2 \circ T_0 '^{40}\circ T_1^{31}(t)=p^{31}t-p^{31}+1~~~~ ,
~~~~\pi_2 \circ  T_1 '^{40}\circ T_0^{31}(t)=p^{31}t+(p^{31}-1)s$\\
of course, if it is needed.

\textbf{Case3.} If $ ~s=\gamma ^{-1} s_{2}~$ or $ ~s=\gamma s_{1}$, then\\
$~~~~~~~~~~min \big\{\frac{p^{i}}{q^{j}} ~\mid~~p^{i}>q^{j} ~\big\}= \frac{p^{7}}{q^{9}}=\gamma~~~ ,~~ ~ max\big\{\frac{p^{i}}{q^{j}}~\mid~~p^{i}<q^{j}~\big\}=\frac{p^{24}}{q^{31}}=\gamma^{-1},$\\
therefore point $(\gamma ^{-1} s_{2},t)\in \mathcal{L} $ passes on line $s=\frac{p^{7}}{q^{9}}\cdot \gamma ^{-1} s_{2}=s_2$ and it had to fall in set ${s_2}\times (-s_2,~1)$ since satisfies  condition of the Lemma 2. By using
 the operators $T_0 '^{40}\circ T_1^{31}$ and  $   T_1 '^{40}\circ T_0^{31}$  and Lemma 1, we transfer  point $(s,t)$ to $ \mathcal{L}^\circ $.

The same  situation is valid for point $ (\gamma s_{1},t)\in\mathcal{L} $.

This completes the proof of the theorem. $\Box$

Replacing  $K$ and $K'$ to the forms  of $C_\alpha$ and $ C_\beta$, respectively and also appropriate selection of the recurrent  set in Theorem 2 yields below proposition.

\begin{proposition}
Set $C_\alpha- \lambda C_\beta$ contains at least one interval for each  $\lambda \in \mathbb{R}^*$.
\end{proposition}
We close this paper by posing a question which inspired by Theorem 2 and Palis conjecture.

\begin{openprob}
Does there exist an open and dense subset in   the mysterious region
$$\Omega:=\Big\{ (C_\alpha, C_\beta)~\Big|~~ HD (C_\alpha)+HD(C_\beta)>1~~and~~\tau(C_\alpha)\cdot\tau(C_\beta)<1~\Big\},$$
such that their elements have stable intersection?
\end{openprob}

\begin{acknowledgment}

I offer my sincerest gratitude to C. G Moreira for posing to me the problem and his valuable comments and suggestions from the very early stage of this research in IMPA.

\end{acknowledgment}

\bigskip
\section*{References}
\small
$~~~~$[HMP] B. Honary, C. G. Moreira, M. Pourbarat,  \textit{Stable intersections of affine Cantor sets},
 Bull Braz Math Soc,
 {\bf 36} (2005), no.~3, 363--378.

[M1] C. G. Moreira, \textit{Stable intersections of Cantor sets and
homoclinic bifurcations}, Ann. Inst. H. Poincar\'e Anal. Non Lin\'eaire
{\bf 13} (1996), no.~6, 741--781.

[M2] C. G. Moreira, \textit{Sums of regular Cantor sets, dynamics and applications to number theory}, Periodica Mathematica Hungarica
{\bf 37} (1998), no.~1-3, 55-63.

[MY] C. G. Moreira\ and\ J.-C. Yoccoz, \textit{Stable intersections of
regular Cantor sets with large Hausdorff dimension}, Ann. of Math.
{\bf 154} (2001), no.~1, 45--96.

[MO] P. Mendes\ and\ F. Oliveira, \textit{On the topological structure of
the arithmetic sum of two Cantor sets},
 Nonlinearity {\bf 7} (1994),no.~2, 329--343.

[N] S. Newhouse, \textit{Non density of Axiom A(a) on $S^2$},
 Proc. A.M.S. Symp. Pure Math {\bf 14} (1970), 191-202.

[P1] J. Palis, \textit{A global view of dynamics and conjecture on
denseness of finitude of attractors}, Ast\'erisque No. 261 (2000),
 xiii--xiv, 335--347.

[P2] J. Palis, \textit{A global perspective for non-conservative dynamics}, Ann. 1.H. Poincar\'e AN {\bf 22} (2005),
, 485-507.

[P3] M. Pourbarat, \textit{On the arithmetic difference of middle Cantor sets}, Submitted.

[PT] J. Palis\ and\ F. Takens, {\it Hyperbolicity and sensitive
 chaotic dynamics at homoclinic bifurcations}, Cambridge Univ. Press,
 Cambridge, 1993.

[S] A. Sannami,  \textit{An example of a regular Cantor
 set whose difference set is a Cantor set with positive measure},
  Hokkaido Math. J. {\bf 21} (1992), no.~1, 7--24.

\end{document}